\documentclass[10pt]{article}
\usepackage{mathrsfs}
\usepackage{amsthm,amsmath,amssymb,anysize}
\newtheorem{lemma}{Lemma}[section]
\newtheorem{theorem}[lemma]{Theorem}
\newtheorem{remark}[lemma]{Remark}
\newtheorem{proposition}[lemma]{Proposition}
\newtheorem{definition}[lemma]{Definition}
\newtheorem{corollary}[lemma]{Corollary}
\usepackage[numbers,sort&compress]{natbib}
\marginsize{3.6cm}{3.6cm}{3.6cm}{3cm}
\numberwithin{equation}{section}

\title{\textsf{On the structures of split $\delta$ Jordan-Lie algebras}}

\author{\textsc{Yan Cao$^{1,2}$\footnote{Supported by Scientific Research Fund of Heilongjiang Provincial Education Department
 (No. 12541184). }, \;Liangyun Chen$^{1}$\footnote{Correspondence:
chenly640@nenu.edu.cn.}
 \footnote{Supported by  NNSF of China (Nos. 11171055 and 11471090),  NSF of  Jilin province (No. 201115006).}
  }\\
  \\
\small $^{1}$School of Mathematics and Statistics, Northeast Normal University \\
 \small Changchun  130024, P. R. China\\
\small $^{2}$Department of Basic Education,  Harbin University of Science and Technology\\
\small Rongcheng Campus, Rongcheng 264300, P. R. China\\
 }
\date{ }

\begin{document}
\maketitle
\begin{quotation}
{\small\noindent\textbf{Abstract:}
We study the structures of arbitrary split $\delta$ Jordan-Lie algebras with
symmetric root systems. We show that any of such algebras $L$ is of the form
$L = U + \sum\limits_{[j] \in \Lambda/\sim}I_{[j]}$ with $U$ a subspace of  $H$ and any $I_{[j]}$, a well described
ideal of $L$, satisfying $[I_{[j]}, I_{[k]}] = 0$ if $[j]\neq [k]$. Under certain conditions, the simplicity of $L$
is characterized and it is shown that $L$ is the direct sum of the family of its minimal ideals,
each one being a simple split $\delta$ Jordan-Lie algebra with a symmetric root system and having all
its nonzero roots connected.}   \\

  \end{quotation}
\setcounter{section}{-1}
\section{Introduction}
\quad $\delta$ Jordan-Lie algebras were  introduced by Susumu Okubo in 1997 \cite{BL4}. The case of $\delta=1$  implies $\delta$ Jordan-Lie algebras
are Lie algebras and the other case of $\delta=-1$ gives  Jordan-Lie algebras. So a question arises whether some known results on Lie algebras can be extended to the framework of $\delta$ Jordan-Lie algebras. $\delta$ Jordan-Lie algebras are the natural generalization of Lie algebras
and have important applications.

In the present paper, we are interested in studying the structures of  arbitrary  $\delta$ Jordan-Lie algebras by focussing on the split ones. The class of the split ones is specially related to addition quantum
numbers, graded contractions, and deformations. For instance, for a physical system which displays
a symmetry of $L$, it is interesting to know in detail the structure of the split decomposition
because its roots can be seen as certain eigenvalues which are the additive quantum numbers
characterizing the state of such system.  Recently, in \cite{BL52, BL52567, BL5234, BL528}, the structures of arbitrary split Lie algebras, arbitrary split
  Leibniz algebras and arbitrary split  Lie triple systems have been determined by the techniques of connections of roots.

Throughout this paper,  $\delta$ Jordan-Lie algebras $L$ are considered of arbitrary dimension and over
an arbitrary field $\mathbb{K}$. In section 1,  we introduce the concept
of split $\delta$ Jordan-Lie algebras as the natural extension of  the split Lie algebras (see \cite{BL52}).
 In section 2,  we improve the techniques of connections of roots introduced for split Lie algebras, so as to develop a theory of connections of roots
for split  $\delta$ Jordan-Lie algebras. Finally, in
 section 3, and under certain conditions, the simplicity of a split $\delta$ Jordan-Lie  algebra $L$ is characterized and it is shown that
$L$ is the direct sum of the family of its minimal ideals, each one being a simple
  split $\delta$ Jordan-Lie  algebra.

\section{Preliminaries}
\begin{definition}{\rm\cite{BL4}}\label{111}
 A $\delta$ Jordan-Lie algebra $L$ is a vector space over a field $\mathbb{K}$ endowed with a bilinear map $[\cdot,\cdot]:L\times L\rightarrow L$ satisfying

$(1)$ $[x,y]=-\delta[y,x]$,

$(2)$ $[x,[y,z]]=\delta[[x,y],z]+\delta[y,[x,z]]$,

\noindent for  $\delta=\pm 1$, $x,y,z\in L$.
 \end{definition}

\begin{remark}{\rm\cite{BL4}}\label{t16789}
A $\delta$ Jordan-Lie algebra  $L$ is called a Lie algebra if $\delta=1$, and a $\delta$ Jordan-Lie algebra  $L$ is called a Jordan-Lie algebra if $\delta=-1$.
 \end{remark}

\begin{definition}{\rm\cite{BL4}}\label{333tb16789}
 Let $L$ be a $\delta$ Jordan-Lie algebra. We define $\rm ad$$:$ $L\rightarrow$ $\rm End$$L$ by $( \mathrm{ad}x)y:= \delta [x,y],$ for  $x,y\in L$.
 \end{definition}

\begin{definition}\label{thmcltb16789345} A splitting Cartan subalgebra $H$ of a  $\delta$ Jordan-Lie algebra $L$ is defined as a maximal abelian
subalgebra $(\rm MASA)$ of $L$ satisfying that the adjoint mappings $\rm ad$$(h)$ for $h \in H$ are simultaneously
diagonalizable. If $L$ contains a splitting Cartan subalgebra $H$, then $L$ is called a split $\delta$ Jordan-Lie
algebra.
\end{definition}\label{thmcltb16789345}

This means that we have a root decomposition $L$ = $H$ $\oplus (\oplus_{\alpha \in \Lambda}L_{\alpha})$  where
$L_{\alpha} =\{v_{\alpha}\in L: [h, v_{\alpha}] = \alpha(h)v_{\alpha} \ for\ any\ h\in H\}$ for a linear function $\alpha \in H^{\ast}
$ and
$\Lambda:= \{\alpha \in H^{\ast}
\setminus \{0\}: L_{\alpha}\neq 0\}$ is the corresponding root system. The subspaces $L_{\alpha}$ for
$\alpha \in H^{\ast}$
 are called root spaces of $L$ $($with respect to $H$$)$ and the elements $\alpha \in \Lambda\cup \{0\}$ are
called roots of $L$ with respect to $H$.

\begin{definition}\label{thmcltb16789345789}
 A root system $\Lambda$ is called symmetric if it satisfies that
$\alpha \in \Lambda$ implies $-\alpha \in \Lambda$.
\end{definition}\label{thmcltb16789345789}

\begin{lemma}\label{thmcltb16789345789345}
If $[L_{\alpha},L_{\beta}]\neq 0$ with $\alpha,\beta \in \Lambda \cup \{0\}$, then $\delta(\alpha+\beta) \in \Lambda\cup \{0\}$ and $[L_{\alpha},L_{\beta}]\subseteq L_{\delta(\alpha+\beta)}$.
\end{lemma}\label{thmcltb16789345789345}

\noindent\textit{Proof.} \quad For any $x\in L_{\alpha}$, $y\in L_{\beta}$ and $h\in H$, by $(2)$ of Definition $\ref{111}$, we have $[h,[x,y]]=\delta][h,x],y]+\delta[x,[h,y]]$=$\delta(\alpha+\beta)(h)[x,y]$.
 Therefore, $[L_{\alpha},L_{\beta}]\subseteq L_{\delta(\alpha+\beta)}$.

\begin{definition}\label{thmcltb16789345789345}
 A subset $\Lambda_{0}$ of $\Lambda$  is called a root subsystem $($relative to $H$$)$ if it is symmetric and $\alpha, \beta$ $\in \Lambda_{0}$,
$\delta(\alpha + \beta) \in  \Lambda$ implies $\delta(\alpha + \beta) \in  \Lambda_{0}$.
\end{definition}\label{thmcltb16789345789345}

Let $\Lambda_{0}$  be a root subsystem of $\Lambda$. We define
 $$H_{\Lambda_{0}} := \mathrm{span}_{\mathbb{K}}\{[L_{\alpha},L_{-\alpha}]:  \alpha \in  \Lambda_{0}\}$$
 and
$V_{\Lambda_{0}} := \oplus_{\alpha \in \Lambda_{0}}L_{\alpha}$. It is straightforward to verify $L_{\Lambda_{0}} := H_{\Lambda_{0}}
\oplus V_{\Lambda_{0}}$ is a  subalgebra
of $L$. We will say that $L_{\Lambda_{0}}$ is a  subalgebra associated to the root subsystem $\Lambda_{0}$.

\section{Decompositions}
\quad In the following, unless otherwise stated,  $L$ denotes a split $\delta$  Jordan-Lie algebra with a symmetric
root system, and  $L$ = $H\oplus (\oplus_{\alpha \in \Lambda}L_{\alpha})$ the corresponding root decomposition.  We begin the study of split $\delta$  Jordan-Lie algebras by developing the concept of connections of roots.

\begin{definition}\label{thmcltb1678934578934511}  Let $\alpha$ and $\beta$ be two nonzero roots, we shall say that $\alpha$ and $\beta$ are connected if there exist $\alpha_{1},\cdots,\alpha_{n}\in \Lambda$  such that

\noindent $\rm (1)$ $ \alpha_{1}=\alpha,$

\noindent $\rm (2)$ $\{\alpha_{1}, \delta(\alpha_{1}+\alpha_{2}), \alpha_{1}+\alpha_{2}+\delta\alpha_{3},\cdots, \delta^{n-2}\alpha_{1}+\sum_{i=1}^{n-2}\delta^{n-1-i}\alpha_{i+1}\}\subset \Lambda,$ \quad $(n\geq 2)$

\noindent $\rm (3)$ $\delta^{n-1}\alpha_{1}+\sum_{i=1}^{n-1}\delta^{n-i}\alpha_{i+1}=\pm\beta$. \quad $(n\geq 1)$

We shall also say
that $\{\alpha_{1},\cdots,\alpha_{n}\}$ is a connection from $\alpha$ to $\beta$.
\end{definition}\label{thmcltb1678934578934511}

We denote by
$$\Lambda_{\alpha} := \{\beta \in \Lambda: \alpha\ and \ \beta \ are \ connected\}.$$
Let us observe that $\{\alpha\}$ is a connection from $\alpha$ to itself and to $-\alpha$. Therefore $\pm \alpha \in  \Lambda_{\alpha}$.

The next result shows that the connection relation is of equivalence.

\begin{proposition}\label{777}
The relation $\sim$ in $\Lambda$, defined by $\alpha \sim \beta$ if and only if $\beta \in  \Lambda_{\alpha}$, is of equivalence.
\end{proposition}
\noindent\textit{Proof.}  $\{\alpha\}$ is a connection from $\alpha$ to itself and therefore $\alpha \sim \alpha$.

 Let us see the symmetric character of $\sim .$ If $\alpha  \sim  \beta$, $\{\alpha_{1},\cdots,\alpha_{n}\}$ is a connection from $\alpha$ to $\beta$, then
 $$\{\delta^{n-1}\alpha_{1}+\sum_{i=1}^{n-1}\delta^{n-i}\alpha_{i+1}, -\delta\alpha_{n}, -\delta\alpha_{n-1},\cdots,-\delta\alpha_{2}\}$$
 is a connection from $\beta$ to $\alpha$ in case $\delta^{n-1}\alpha_{1}+\sum_{i=1}^{n-1}\delta^{n-i}\alpha_{i+1}=\beta$, and
  $$\{-\delta^{n-1}\alpha_{1}-\sum_{i=1}^{n-1}\delta^{n-i}\alpha_{i+1}, \delta\alpha_{n}, \delta\alpha_{n-1},\cdots,\delta\alpha_{2}\}$$
 is a connection from $\beta$ to $\alpha$ in case $\delta^{n-1}\alpha_{1}+\sum_{i=1}^{n-1}\delta^{n-i}\alpha_{i+1}=-\beta$. Therefore $\beta \sim \alpha$.

Finally, suppose $\alpha \sim \beta$ and $\beta \sim \gamma$, and write $\{\alpha_{1},\cdots,\alpha_{n}\}$ for a connection from $\alpha$ to $\beta$
and $\{\beta_{1},\cdots,\beta_{m}\}$ for a connection from $\beta$ to $\gamma$.  If $m >1$, then
$$\{\alpha_{1},\cdots,\alpha_{n},\beta_{2},\cdots,\beta_{m}\}$$
is a connection from $\alpha$ to $\gamma$ in case $\delta^{n-1}\alpha_{1}+\sum_{i=1}^{n-1}\delta^{n-i}\alpha_{i+1}=\beta$, and $$\{\alpha_{1},\cdots,\alpha_{n},-\beta_{2},\cdots,-\beta_{m}\}$$ in case $\delta^{n-1}\alpha_{1}+\sum_{i=1}^{n-1}\delta^{n-i}\alpha_{i+1}=-\beta.$
If $m = 1$, then $\gamma \in \pm\beta$ and so $\{\alpha_{1},\cdots,\alpha_{n}\}$ is a connection
from $\alpha$ to $\gamma$.  Therefore $\alpha \sim \gamma$. So, $\sim$ is of equivalence.  \hfill$\Box$

\begin{proposition}\label{888}
Let $\alpha$ be a nonzero root. Then the following assertions hold.

\noindent $\rm (1)$ $\Lambda_{\alpha}$ is a root subsystem.

\noindent $\rm (2)$ If $\gamma$ is a nonzero root such that $\gamma \not \in$ $\Lambda_{\alpha}$, then $[L_{\beta},L_{\gamma}]= 0$  and $\gamma ([L_{\beta},L_{-\beta}])= 0$  for
any $\beta \in \Lambda_{\alpha} $.
\end{proposition}

\noindent\textit{Proof.}  $(1)$ If $\beta \in \Lambda_{\alpha}$ then there exists a connection $\{\alpha_{1},\cdots,\alpha_{n}\}$ from $\alpha$ to $\beta$. By Definition 2.1, it is clear that
$\{\alpha_{1},\cdots,\alpha_{n}\}$ also connects $\alpha$ to $-\beta$ and therefore $-\beta \in \Lambda_{\alpha}$. If $\beta, \gamma \in \Lambda_{\alpha}$ and $\delta(\beta+\gamma) \in \Lambda$,
then there exists a connection $\{\alpha_{1},\cdots,\alpha_{n}\}$ from $\alpha$ to $\beta$.
Hence, $\{\alpha_{1},\cdots,\alpha_{n},\gamma\}$ is a
connection from $\alpha$ to $\beta+\gamma$ in case $\delta^{n-1}\alpha_{1}+\sum_{i=1}^{n-1}\delta^{n-i}\alpha_{i+1}=\beta$ and  $\{\alpha_{1},\cdots,\alpha_{n},-\gamma\}$ in case $\delta^{n-1}\alpha_{1}+\sum_{i=1}^{n-1}\delta^{n-i}\alpha_{i+1}=-\beta$. So $\delta(\beta+\gamma) \in  \Lambda_{\alpha}$.

 $(2)$ Let us suppose that there exists $\beta$ $\in \Lambda_{\alpha}$ such that $[L_{\beta},L_{\gamma}]\neq 0$  with $\gamma \not \in \Lambda_{\alpha}$. Then
$\delta(\beta+\gamma) \in \Lambda$  and we have as in Proposition \ref{888} $(1)$ that $\alpha$ is connected to $\delta(\beta+\gamma)$. Since $\Lambda_{\alpha}$ is a root subsystem,
one gets $\gamma \in \Lambda_{\alpha}$, a contradiction. Therefore $[L_{\beta},L_{\gamma}]= 0 $ for any $\beta \in \Lambda_{\alpha}$ and $\gamma \not \in \Lambda_{\alpha}$. As
$-\beta \in \Lambda_{\alpha}$ for any $\beta \in \Lambda_{\alpha}$, we also have that $[L_{-\beta},L_{\gamma}]= 0$. By   $(2)$ of Definition $\ref{111}$, one gets $[[L_{\beta},L_{-\beta}],L_{\gamma}]\subset \delta[L_{\beta},[L_{-\beta},L_{\gamma}]]+[L_{-\beta}, [L_{\beta},L_{\gamma}]].$  Therefore  $[[L_{\beta},L_{-\beta}],L_{\gamma}]=0$. So $\gamma ([L_{\beta},L_{-\beta}])= 0$ for
any $\beta \in \Lambda_{\alpha} $ and $\gamma \not \in \Lambda_{\alpha}$.  \hfill$\Box$

\begin{theorem}\label{669}
The following assertions hold.

\noindent $\rm(1)$ For any $\alpha \in \Lambda$, the  subalgebra
$L_{\Lambda_{\alpha}}
= H_{\Lambda_{\alpha}}
\oplus V_{\Lambda_{\alpha}}$
of $L$ associated to the root subsystem $\Lambda_{\alpha}$ is an ideal of $L$.

\noindent $\rm (2)$ If $L$ is simple, then there exists a connection from $\alpha$ to $\beta$ for any $\alpha, \beta \in \Lambda$.
\end{theorem}
\noindent\textit{Proof.}
(1)  Recall that $$H_{\Lambda_{\alpha}} := \mathrm{span}_{\mathbb{K}}\{[L_{\beta},L_{-\beta}]:  \beta \in  \Lambda_{\alpha}\}\subset H$$
 and
$V_{\Lambda_{\alpha}} := \oplus_{\gamma \in \Lambda_{\alpha}}L_{\gamma}.$ In order to complete the proof, it is sufficient to show that
$$[L_{\Lambda_{\alpha}},L]\subset L_{\Lambda_{\alpha}}.$$
It is easy to see that
\begin{align*}
[L_{\Lambda_{\alpha}},L]=&[ \sum\limits_{\beta \in \Lambda_{\alpha}}[L_{\beta},L_{-\beta}]\oplus(\oplus_{\beta \in \Lambda_{\alpha}}L_{\beta}), H\oplus(\oplus_{\gamma \in \Lambda_{\alpha}}L_{\gamma})\oplus(\oplus_{\gamma \not \in \Lambda_{\alpha}}L_{\gamma})]\\
=&[ \sum\limits_{\beta \in \Lambda_{\alpha}}[L_{\beta},L_{-\beta}], H]+ [ \sum\limits_{\beta \in \Lambda_{\alpha}}[L_{\beta},L_{-\beta}], \oplus_{\gamma \in \Lambda_{\alpha}}L_{\gamma}]+[ \sum\limits_{\beta \in \Lambda_{\alpha}}[L_{\beta},L_{-\beta}], \oplus_{\gamma \not \in \Lambda_{\alpha}}L_{\gamma}]\\
+&[\oplus_{\beta \in \Lambda_{\alpha}}L_{\beta}, H]+[\oplus_{\beta \in \Lambda_{\alpha}}L_{\beta}, \oplus_{\gamma \in \Lambda_{\alpha}}L_{\gamma}]+[\oplus_{\beta \in \Lambda_{\alpha}}L_{\beta}, \oplus_{\gamma \not \in \Lambda_{\alpha}}L_{\gamma}].
\end{align*}
Here, it is clear $[ \sum\limits_{\beta \in \Lambda_{\alpha}}[L_{\beta},L_{-\beta}], H]\subset [H, H]=0$. By Proposition \ref{888} (2), one gets
$[ \sum\limits_{\beta \in \Lambda_{\alpha}}[L_{\beta},L_{-\beta}], \oplus_{\gamma \not \in \Lambda_{\alpha}}L_{\gamma}]=0$ and $[\oplus_{\beta \in \Lambda_{\alpha}}L_{\beta}, \oplus_{\gamma \not \in \Lambda_{\alpha}}L_{\gamma}]=0$. Next, we easily obtain $$[ \sum\limits_{\beta \in \Lambda_{\alpha}}[L_{\beta},L_{-\beta}], \oplus_{\gamma \in \Lambda_{\alpha}}L_{\gamma}]\subset [H, \oplus_{\gamma \in \Lambda_{\alpha}}L_{\gamma}]\subset \oplus_{\gamma \in \Lambda_{\alpha}}L_{\gamma}=V_{\Lambda_{\alpha}}$$
and
$$[\oplus_{\beta \in \Lambda_{\alpha}}L_{\beta}, H]\subset \oplus_{\beta \in \Lambda_{\alpha}}L_{\beta}=V_{\Lambda_{\alpha}}.$$
At last, taking into account $[\oplus_{\beta \in \Lambda_{\alpha}}L_{\beta}, \oplus_{\gamma \in \Lambda_{\alpha}}L_{\gamma}]$, we treat two cases.

Case 1:  $\beta+\gamma=0$ for $\beta \in \Lambda_{\alpha}, \gamma \in \Lambda_{\alpha}$. One gets  $[\oplus_{\beta \in \Lambda_{\alpha}}L_{\beta}, \oplus_{\gamma \in \Lambda_{\alpha}}L_{\gamma}]\subset H_{\Lambda_{\alpha}}.$

Case 2: $\beta+\gamma\neq 0$ for $\beta \in \Lambda_{\alpha}, \gamma \in \Lambda_{\alpha}$. Since $\Lambda_{\alpha}$ is a root subsystem, one gets if $[\oplus_{\beta \in \Lambda_{\alpha}}L_{\beta}, \oplus_{\gamma \in \Lambda_{\alpha}}L_{\gamma}]\neq 0$ then $[\oplus_{\beta \in \Lambda_{\alpha}}L_{\beta}, \oplus_{\gamma \in \Lambda_{\alpha}}L_{\gamma}]\subset \oplus_{\delta(\beta+\gamma) \in \Lambda_{\alpha}}L_{\delta(\beta+\gamma)}= V_{\Lambda_{\alpha}} .$

Therefore, $[L_{\Lambda_{\alpha}},L]\subset L_{\Lambda_{\alpha}}$ is proved.

(2) The simplicity of $L$ implies $L_{\Lambda_{\alpha}}
= L$ and therefore $ \Lambda_{\alpha} = \Lambda$.  \hfill$\Box$

\begin{theorem}\label{66777}
For a vector space complement $U$ of $span_{\mathbb{K}}\{{[L_{\alpha},L_{-\alpha}]: \alpha \in \Lambda}\}$ in H, we have
$$L = U + \sum\limits_{[\alpha] \in \Lambda/\sim}I_{[\alpha]},$$
where any $I_{[\alpha]}$ is one of the ideas $L_{\Lambda_{\alpha}}$ of $L$ described in Theorem \ref{669}, satisfying
$[I_{[\alpha]},I_{[\beta]}]=0,$ whenever $[\alpha] \neq [\beta].$
\end{theorem}
\noindent\textit{Proof.} By Proposition \ref{777}, we can consider the quotient set $\Lambda/\sim:= \{[\alpha]: \alpha \in \Lambda\}$.
Let us denote by  $I_{[\alpha]} :=L_{\Lambda_{\alpha}}.$  $I_{[\alpha]}$ is well-defined and is an ideal of $L$ by Theorem \ref{669} (2). Therefore
$$L = U + \sum\limits_{[\alpha] \in \Lambda/\sim}I_{[\alpha]}.$$

Now, given $[\alpha]\neq [\beta]$, the assertion $[I_{[\alpha]},I_{[\beta]}]=0$, follows from Proposition \ref{888} (2).\hfill$\Box$

\begin{definition}
The center of a $\delta$ Jordan-Lie algebra is the set $\mathrm{Z}(L)=\{\upsilon \in L: [\upsilon, L]=0\}$.
\end{definition}

\begin{corollary}\label{667}
If $\mathrm{Z}(L) = 0$ and $[L,L] = L$, then $L$ is the direct sum of the ideals given in Theorem \ref{66777},
$$L =\oplus_{[\alpha] \in \Lambda/\sim}I_{[\alpha]}.$$
\end{corollary}

\noindent\textit{Proof.} From $[L,L] = L$, it is clear that $L = \sum\limits_{[\alpha] \in \Lambda/\sim}I_{[\alpha]}.$
Next, we will show the direct character of the sum. Given $x \in I_{[\alpha]}\cap  \sum\limits_{[\beta] \in \Lambda/\sim \atop \beta \not \sim \alpha}I_{[\beta]},$  using again the equation $[I_{[\alpha]},I_{[\beta]}]=0$ for $[\alpha] \neq [\beta],$ one gets
$[x, I_{[\alpha]} \oplus  \sum\limits_{[\beta] \in \Lambda/\sim \atop \beta \not \sim \alpha}I_{[\beta]} ]=0$, that is, $x \in \mathrm{Z}(L)=0$. Thus $x=0$, as desired.  \hfill$\Box$

\section{The simple components}

\quad In this section we study if any of the components in the decomposition given in Corollary \ref{667}
is simple. Under certain conditions we give an affirmative answer. From now on $char(\mathbb{K})=0$.

\begin{lemma}\label{667188}
Let $L$ be a split $\delta$ Jordan-Lie algebra. For any $\alpha$, $\beta$ $\in \Lambda$  with $\alpha \neq k\beta$, $k \in \mathbb{K}$, there
exists $h_{\alpha,\beta} \in H$ such that $\alpha(h_{\alpha,\beta}) \neq 0$ and $\beta( h_{\alpha,\beta}) = 0$.
\end{lemma}

\noindent\textit{Proof.} This can be proved completely analogously to \cite[Lemma 3.1]{BL52}. \hfill$\Box$

\begin{lemma}\label{66}
Let $L$  be a split  $\delta$ Jordan-Lie algebra. If I is an ideal of $L$ and $x = h_{0}+\sum_{j=1}^{m}e_{\beta_{j}}\in I$, for
$h_{0}\in H$, $e_{\beta_{j}}\in L_{\beta_{j}}$ and $\beta_{j} \neq \beta_{k}$ if $j\neq k$. Then any $e_{\beta_{j}}\in I.$
\end{lemma}

\noindent\textit{Proof.} Let us fix $\beta_{1}$. If $e_{\beta_{1}}=0$ then $e_{\beta_{1}} \in I$. Suppose  $e_{\beta_{1}}\neq 0$.
For any $\beta_{k_{r}}\neq p_{\beta_{1}}, p\in \mathbb{K}$ and $k_{r}\in \{2,\cdots,m\}$, Lemma \ref{667188} gives us $h_{\beta_{1},\beta_{k_{r}}} \in H$
satisfying     $\beta_{1}(h_{\beta_{1},\beta_{k_{r}}})\neq 0$ and $\beta_{k_{r}}(h_{\beta_{1},\beta_{k_{r}}})=0$.  From  here,

\begin{equation}\label{11}
[[\cdots[[x, h_{\beta_{1},\beta_{k_{2}}}],h_{\beta_{1},\beta_{k_{3}}}],\cdots],h_{\beta_{1},\beta_{k_{s}}}]=\delta(p_{1}e_{\beta_{1}}+\sum_{t=1}^{u}p_{k_{t}}e_{k_{t}\beta_{1}})\in I,
\end{equation}

\noindent where $p_{1}$, $k_{t}\in \mathbb{K}-\{0\},k_{t}\neq 1$ and $p_{k_{t}} \in \mathbb{K}$.

If any $p_{k_{t}}=0$, $t=1,\cdots,u,$ then $p_{1}e_{\beta_{1}}\in I$ and so $e_{\beta_{1}}\in I$. Let us suppose some $p_{k_{t}}\neq 0$ and write ($\ref{11}$) as
\begin{equation}\label{12}
p_{1}e_{\beta_{1}}+\sum_{t=1}^{v}p_{k_{t}}e_{k_{t}\beta_{1}}\in I,
\end{equation}
where $p_{1}, k_{t},p_{k_{t}}\in \mathbb{K}-\{0\}$, $k_{t} \neq 1, v\leq u$.

Let $h \in H$  such that $\beta_{1}(h)\neq 0$. Then

$$[h,  p_{1}e_{\beta_{1}}+\sum_{t=1}^{v}p_{k_{t}}e_{k_{t}\beta_{1}}]=p_{1}\beta_{1}(h)e_{\beta_{1}}+\sum_{t=1}^{v}p_{k_{t}}k_{t}\beta_{1}(h)e_{k_{t}\beta_{1}}\in I,$$

\noindent and so

\begin{equation}\label{13}
p_{1}e_{\beta_{1}}+\sum_{t=1}^{v}p_{k_{t}}k_{t}e_{k_{t}\beta_{1}}\in I,\quad k_{t}\neq 1.
\end{equation}

From ($\ref{12}$) and ($\ref{13}$), it follows easily that

\begin{equation}\label{14}
q_{1}e_{\beta_{1}}+\sum_{t=1}^{w}q_{k_{t}}e_{q_{t}\beta_{1}}\in I,
\end{equation}

\noindent where $q_{1},q_{k_{t}} \in \mathbb{K}-\{0\}$, $q_{t} \in \{k_{t}: t=1,\cdots,v\}$ and $w<v$.

Following this process $($multiply ($\ref{14}$) with $h$ and compare the result with ($\ref{14}$) taking into account $q_{t}\neq 1$, and so
on$)$, we obtain $e_{\beta_{1}}\in I$.
 The same argument holds for any $\beta_{j},j\neq 1.$  \hfill$\Box$

\begin{definition}\label{thmcltb167893457893451156}
We say that a split $\delta$ Jordan-Lie algebra $L$ is root-multiplicative if $\alpha, \beta, \alpha+\beta \in \Lambda$ implies $[L_{\alpha},L_{\beta}]\neq 0.$
\end{definition}

\begin{lemma}\label{567}
Let $L$ be a split $\delta$ Jordan-Lie algebra with $\mathrm{Z}(L) = 0$ and $I$ be an ideal of $L$. If $I\subset H$ then $I=0$.
\end{lemma}
\noindent\textit{Proof.} Suppose there exists a nonzero ideal of $L$ such that $I\subset H$.  On the one hand, $$[I,  \oplus_{\alpha \in \Lambda}L_{\alpha}]\subset I \subset H.$$
 On the other hand,
 $$[I, \oplus_{\alpha \in \Lambda}L_{\alpha}]\subset [H, \oplus_{\alpha \in \Lambda}L_{\alpha}]\subset \oplus_{\alpha \in \Lambda}L_{\alpha}.$$  So,
   $$[I, \oplus_{\alpha \in \Lambda}L_{\alpha}]\subset H\cap  \oplus_{\alpha \in \Lambda}L_{\alpha}=0. $$
   From $[I,H]=0$ and  $[I, \oplus_{\alpha \in \Lambda}L_{\alpha}]=0$, one gets
$$[I,H\oplus(\oplus_{\alpha \in \Lambda}L_{\alpha})]= [I, L] =0.$$
 Thus $I\subset \mathrm{Z}(L)=0$, a contradiction.\hfill$\Box$

\begin{theorem}\label{9}
Let $L$ be a split $\delta$ Jordan-Lie algebra.  If $L$ is   root-multiplicative, with $\mathrm{Z}(L) = 0$, $[L,L] = L$ and satisfying
$\mathrm{dim}L_{\alpha} = 1$ for any $\alpha \in \Lambda$. Then $L$ is simple if and only if it has all its nonzero roots
connected.
\end{theorem}

\noindent\textit{Proof.} The first implication is Theorem \ref{669} (2). Let us prove the converse.  Let $I$ be a
nonzero ideal of $L$. In order to complete the proof, it is sufficient to show $L\subset I$. By $[L,L] = L$, one gets
$$L=\sum\limits_{\alpha \in \Lambda}[L_{\alpha},L_{-\alpha}]\oplus(\oplus_{\beta \in \Lambda}L_{\beta}).$$ Next, we  check
$$\sum\limits_{\alpha \in \Lambda}[L_{\alpha},L_{-\alpha}]\oplus(\oplus_{\beta \in \Lambda}L_{\beta})\subset I.$$
By Lemma \ref{567}, we can find $0 \neq x \in  I$  such that $$x = h_{0} + \sum_{j=1}^{m}e_{\beta_{j}} \in I,$$
 where
$h_{0} \in H$, $e_{\beta_{j}}
\in  L_{\beta_{j}}$, $\beta_{j} \neq \beta_{k}$  if $j \neq k$ and satisfying some $e_{\beta_{j}}
\neq 0$.

Let us choose such an
$x \in I,$  and fix any $\beta_{j_{0}}$, $j_{0} \in \{1,\cdots , m\}$ such that $e_{\beta_{j_{0}}}\neq 0.$
By applying Lemma \ref{66},
$e_{\beta_{j_{0}}} \in I$
   and as dim$L_{\beta_{j_{0}}}=1,$
 one gets $L_{\beta_{j_{0}}}\subset I$. We conclude  that $L_{-\beta_{j_{0}}}\subset I.$ Indeed, since $\beta_{j_{0}}\neq 0$ and $[L,L]=L$, there exists $[e_{\gamma} , e_{-\gamma}] \neq 0$, $e_{\pm\gamma} \in L_{\pm\gamma}$,
$\gamma \in \Lambda$, such that $\beta_{j_{0}}([e_{\gamma}, e_{-\gamma}]) \neq 0$. In the following, we treat two cases.

Case 1:  $\gamma \in \pm\beta_{j_{0}}$. As $0 \neq [e_{\beta_{j_{0}}},
 e_{-\beta_{j_{0}}}] \in I$, one has
$$ e_{-\beta_{j_{0}}}
=  -\beta_{j_{0}} ([e_{\beta_{j_{0}}},
 e_{-\beta_{j_{0}}}])^{-1}[[e_{\beta_{j_{0}}},
 e_{-\beta_{j_{0}}}], e_{-\beta_{j_{0}}}] \in I,$$
  and so $$L_{-\beta_{j_{0}}}\in I.$$

 Case 2: $\gamma \not \in \pm\beta_{j_{0}}$. As $\beta_{j_{0}}$ and $\gamma$
 are connected, the root-multiplicativity of $L$ and the assumption dim$L_{\alpha} = 1$
for any $\alpha \in \Lambda$, one gets a connection $\{\alpha_{1}, \cdots , \alpha_{r}\}$ from $\beta_{j_{0}}$ to $\gamma$ such that $$\alpha_{1}=\beta_{j_{0}},
\delta(\alpha_{1}+\alpha_{2}),\cdots, \delta^{r-1}\alpha_{1}+\sum_{i=1}^{r-1}\delta^{r-i}\alpha_{i+1}=\pm\gamma \in \Lambda,$$
with $$[L_{\alpha_{1}},L_{\alpha_{2}}]=L_{\delta(\alpha_{1}+\alpha_{2})}, $$ $$[[L_{\alpha_{1}},L_{\alpha_{2}}],L_{\alpha_{3}}]=L_{\alpha_{1}+\alpha_{2}+\delta\alpha_{3}},\cdots,[[\cdots[[L_{\alpha_{1}},L_{\alpha_{2}}],L_{\alpha_{3}}],\cdots,],L_{\alpha_{r}}]=L_{\epsilon\gamma},$$
where $\epsilon \in \pm1$, and deduce that either $L_{\gamma}\subset I$ or $L_{-\gamma}\subset I.$
In both cases $[L_{\gamma} ,L_{-\gamma}] \subset I$ and so $[e_{\gamma}, e_{-\gamma}]\in  I$. As
given any $e_{-\beta_{j_{0}}} \in L_{-\beta_{j_{0}}}$,
 one has $$e_{-\beta_{j_{0}}}
=  -\beta_{j_{0}} ([e_{\gamma},
 e_{-\gamma}])^{-1}[[e_{\gamma},
 e_{-\gamma}], e_{-\beta_{j_{0}}}] \in I,$$
 then we conclude $$L_{-\beta_{j_{0}}}\subset I,$$
 and so $$[L_{\beta_{j_{0}}}
,L_{-\beta_{j_{0}}}
]\subset I.$$  Given any $\tau \in \Lambda$, $\tau \neq \pm_{\beta_{j_{0}}} $, arguing
as above, we have some $L_{\epsilon\tau}\subset I, \epsilon \in \pm1$, and from here $L_{-\epsilon\tau}\subset I$. It is clear that
$[L_{\tau},L_{-\tau}] \subset I$. This proves
$\sum\limits_{\alpha \in \Lambda}[L_{\alpha},L_{-\alpha}]\oplus(\oplus_{\beta \in \Lambda}L_{\beta})\subset I$. So we conclude that $L$ is simple. \hfill$\Box$

\begin{theorem}\label{66712}
 Let $L$ be a split $\delta$ Jordan-Lie algebra. If $L$ be  root-multiplicative, with $\mathrm{Z}(L) = 0$, $[L,L] = L$ and satisfying
$\mathrm{dim}L_{\alpha} = 1$ for any $\alpha \in \Lambda$. Then $L$
is the direct sum of the family of its minimal ideals,
each one being a simple split $\delta$ Jordan-Lie algebra with a symmetric root system and having all its
nonzero roots connected.
\end{theorem}

\noindent\textit{Proof.} By Corollary \ref{667}, $L =\oplus_{[\alpha_{0}] \in \Lambda/\sim}I_{[\alpha_{0}]}$
 is the direct sum of the ideals $I_{[\alpha_{0}]}$=$H_{\Lambda_{\alpha_{0}}}\oplus V_{\Lambda_{\alpha_{0}}}$=$(\sum_{\alpha \in [\alpha_{0}]}[L_{\alpha},L_{-\alpha}])$
 $\oplus(\oplus_{\alpha\in [\alpha_{0}]}L_{\alpha})$ having any $I_{[\alpha_{0}]}$ its root
   system $\Lambda_{\alpha_{0}}$, with all of its roots connected.
Taking into account the facts $\Lambda_{\alpha_{0}}=-\Lambda_{\alpha_{0}}$ and $[L_{\Lambda_{\alpha_{0}}},L_{\Lambda_{\alpha_{0}}}]\subset L_{\Lambda_{\alpha_{0}}}$
$($see Proposition \ref{888} (1)$)$, we easily
deduce that $\Lambda_{\alpha_{0}}$ has all of its roots $\Lambda_{\alpha_{0}}$-connected $($connected through roots in $\Lambda_{\alpha_{0}})$. We also have
that any of the $I_{[\alpha_{0}]}$ is root-multiplicative as consequence of the root-multiplicativity of $L$. Clearly
dim$(I_{[\alpha_{0}]})_{\alpha} =1, \alpha \in [\alpha_{0}]$, and finally $\mathrm{Z}_{I_{[\alpha_{0}]}}I_{[\alpha_{0}]}=0$ $($ where  $\mathrm{Z}_{I_{[\alpha_{0}]}}I_{[\alpha_{0}]}$denotes the center of $I_{[\alpha_{0}]}$ in
 $I_{[\alpha_{0}]}$$)$, as a consequence of  $[I_{[\alpha_{0}]}, I_{[\beta_{0}]}]=0$ if $[\alpha_{0}] \neq [\beta_{0}]$ $($Theorem \ref{66777}$)$ and $\mathrm{Z}(L)=0$. We can apply Theorem $\ref{9}$ to any $I_{[\alpha_{0}]}$ so as to conclude $I_{[\alpha_{0}]}$ is  a simple $\delta$ Jordan-Lie algebra.

 Next, similar to \cite[Theorem 3.5]{BL52}, it can easily show that any of the simple $\delta$ Jordan-Lie algebras $($see Corollary \ref{667}$)$,  $I_{[\alpha_{0}]}$=$(\sum\limits_{\alpha \in [\alpha_{0}]}[L_{\alpha},L_{-\alpha}])$
 $\oplus(\oplus_{\alpha \in [\alpha_{0}]}L_{\alpha})$ is split.  \hfill$\Box$

\vspace{0.5cm}

\end{document}